\documentclass{amsart}
\usepackage{amssymb,amsmath,amsthm,mathrsfs,graphics,hyperref}
\newcommand{\includefigure}[3]{{
  \begin{center}
  \resizebox{#1}{#2}{\includegraphics{{#3}}}
  \end{center}}}

\newtheorem{theorem}{Theorem}
\newtheorem{proposition}[theorem]{Proposition}
\newtheorem{lemma}[theorem]{Lemma}

\theoremstyle{definition}
\newtheorem{definition}{Definition}

\newcommand{\MGW}{\mathcal{G}}
\newcommand{\MRW}{\mathcal{R}}
\newcommand{\DD}{\mathcal{D}}
\newcommand{\Cc}{\mathbb{C}}
\newcommand{\Nn}{\mathbb{N}}

\newcommand{\0}{\emptyset}
\newcommand{\commentout}[1]{}
\newcommand{\sm}{\setminus}
\newcommand{\st}{~|~}
\newcommand{\Sym}{\mathfrak{S}}

\begin{document}

\title{Updown numbers and the initial monomials of the slope variety}
\author{Jeremy L.\ Martin and Jennifer D.\ Wagner}

\address{
Department of Mathematics,
University of Kansas,
405~Snow Hall,
1460~Jayhawk Blvd.,
Lawrence,~KS 66045, USA}
\email{jmartin@math.ku.edu}

\address{
Department of~Mathematics and~Statistics,
Washburn~University,
1700~SW~College~Ave.,
Topeka,~KS 66621, USA}
\email{jennifer.wagner1@washburn.edu}


\begin{abstract}
Let $I_n$ be the ideal of all algebraic relations on
the slopes of the $\binom{n}{2}$ lines formed by placing $n$ points in a plane
and connecting each pair of points with a line.  Under each of two natural term
orders, the initial ideal of $I_n$ is generated by monomials corresponding
to permutations satisfying a certain pattern-avoidance condition.  We show bijectively that these
permutations are enumerated by the updown (or Euler) numbers,
thereby obtaining a formula for the number of generators of
the initial ideal of $I_n$ in each degree.
\end{abstract}
\date{May 28, 2009}
\thanks{First author partially supported by an NSA Young Investigator's Grant}
\maketitle


The symbol $\Nn$ will denote the set of positive integers.
For integers $m\leq n$, we put $[n]=\{1,2,\dots,n\}$ and
$[m,n]=\{m,m+1,\dots,n\}$.  The set of all permutations of an integer set
$P$ will be denoted $\Sym_P$, and the $n^{th}$ symmetric group is $\Sym_n$ ($=\Sym_{[n]}$).
We will write each permutation $w\in\Sym_P$
as a word with $n=|P|$ digits, $w=w_1\dots w_n$, where
$\{w_1,\dots,w_n\}=P$.  The symbol $w_i^{-1}$ denotes the position
of digit $i$ in $w$; that is, $w_i^{-1}=j$ if and only if $w_j=i$.
If necessary for clarity, we will separate the digits
with commas.  Concatenation will also be denoted with commas; for instance,
if $w=12$ and $w'=34$, then $(w,w',5)=12345$.  The \emph{reversal} $w^*$ of $w_1w_2\dots w_{n-1}w_n$
is the word $w_nw_{n-1}\dots w_2w_1$.  A \emph{subword} of a permutation $w\in\Sym_P$ is a
word $w[i,j]=w_iw_{i+1}\cdots w_j$,
where $[i,j]\subseteq[n]$.  The subword is \emph{proper} if $w[i,j]\neq w$.
We write $w\approx w'$ if the digits of $w$ are in the same relative order as those of $w'$;
for instance, $58462\approx 35241$.

\begin{definition} 
\label{G-word}
Let $P\subset\Nn$ with $n=|P|\geq 2$.
A permutation $w\in\Sym_P$ is a \emph{G-word} if
it satisfies the two conditions
  \begin{itemize}
  \item[\textbf{(G1)}] $w_1 = \max(P)$ and $w_n = \max(P\sm\{w_1\})$; and
  \item[\textbf{(G2)}] If $n\geq 4$, then $w_2>w_{n-1}$.
  \end{itemize}
It is an \emph{R-word} if it satisfies the two conditions
  \begin{itemize}
  \item[\textbf{(R1)}] $w_1 = \max(P)$ and $w_n = \max(P\sm\{w_1\})$; and
  \item[\textbf{(R2)}] If $n\geq 4$, then $w_2<w_{n-1}$.
  \end{itemize}
A G-word (resp., an R-word) is \emph{primitive} if for every proper subword
$x$ of length~$\geq 4$, neither $x$ nor $x^*$ is a G-word (resp., an R-word).
The set of all primitive G-words (resp., on $P\subset\Nn$, or on $[n]$)
is denoted $\MGW$ (resp., $\MGW_P$, or $\MGW_n$).  The sets
$\MRW$, $\MRW_P$, $\MRW_n$ are defined similarly.
\end{definition}

For example, the word $53124$ is a G-word, but not a primitive one, because
it contains the G-word $3124=(4213)^*$ as a subword.
The primitive G- and R-words of lengths up to 6 are as follows:

\begin{equation}
\label{small-examples}
\begin{aligned}
\MGW_2 &= \{21\},\\
\MGW_3 &= \{312\},\\
\MGW_4 &= \{4213\},\\
\MGW_5 &= \{52314,\ 53214\},\\
\MGW_6 &= \{623415,\ 624315,\ 642315,\ 634215,\ 643215\},\\
\MRW_2 &= \{21\},\\
\MRW_3 &= \{312\},\\
\MRW_4 &= \{4123\},\\
\MRW_5 &= \{51324,\ 52134\},\\
\MRW_6 &= \{614235,\ 624135,\ 623145,\ 621435,\ 631245\}.
\end{aligned}
\end{equation}

Clearly, if $w\approx w'$, 
then either both $w$ and $w'$ are (primitive) G- (R-)words, or neither are;
therefore, for all $P\subset\Nn$, the set $\MGW_P$ is determined by
(and in bijection with) $\MGW_{|P|}$.

These families of permutations arose in \cite{Slopes} in the following way.
Let $p_1=(x_1,y_1),\dots,p_n=(x_n,y_n)$ be points in $\Cc^2$ with distinct
$x$-coordinates, let $\ell_{ij}$ be the unique line through $p_i$ and $p_j$,
and let $m_{ij}=(y_j-y_i)/(x_j-x_i)\in\Cc$ be the slope of $\ell_{ij}$.
Let $A=\Cc[m_{ij}]$, and let $I_n\subset A$ be the
ideal of algebraic relations on the slopes $m_{ij}$ that hold for all choices of the points $p_i$.
Order the variables of $A$ lexicographically by their subscripts: $m_{12}<m_{13}<\cdots<m_{1n}<m_{23}<\cdots$.
Then \cite[Theorem~4.3]{Slopes}, with respect to graded lexicographic order
on the monomials of $A$, the initial ideal of $I_n$
is generated by the squarefree monomials $m_{w_1,w_2} m_{w_2w_3}\cdots
m_{w_{r-1}w_r}$, where $\{w_1,\dots,w_r\}\subseteq[n]$, $r\geq 4$, and
$w=w_1w_2\cdots w_r$ is a primitive G-word.
Consequently, the number of degree-$d$ generators of the initial ideal of $I_n$ is
  \begin{equation}
  \label{num-gen}
  \binom{n}{d+1}|\MGW_{d+1}|.
  \end{equation}
Similarly, under reverse lex order (rather than graded lex order) on~$A$,
the initial ideal of~$I_n$ is generated by the squarefree monomials
corresponding to primitive R-words.

It was noted in \cite[p.~134]{Slopes} that the first several values of the
sequence $|\MGW_3|,|\MGW_4|,\dots$
coincide with the \emph{updown numbers} (or \emph{Euler numbers}):
  $$1,1,2,5,16,61,272,\dots.$$
This is sequence A000111 in
the Online Encyclopedia of Integer Sequences~\cite{EIS}.
The updown numbers enumerate (among other things) 
the \emph{decreasing 012-trees} \cite{Callan,Donaghey},
which we now define.


\begin{definition}
\label{decreasing-tree}
A \emph{decreasing 012-tree} is a rooted tree, with vertices labeled by distinct
positive integers, such that (i) every vertex has either 0, 1, or 2 children;
and (ii) $x<y$ whenever $x$ is a descendant of $y$.  The set of all decreasing
012-trees with vertex set $P$ will be denoted $\DD_P$.
We will represent rooted trees by the recursive notation $T=[v,T_1,\dots,T_n]$,
where the $T_i$ are the subtrees rooted at the children of $v$.
Note that reordering the $T_i$ in this notation does not change the tree~$T$.
For instance, $[6,[5,[4],[2]],[3,[1]]]$ represents the decreasing 012-tree shown below.
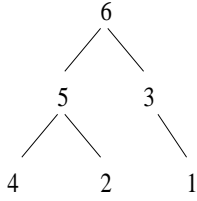
\includefigure{1.0in}{1.0in}{d012tree}
\end{definition}

The purpose of this note is to verify that the updown numbers do indeed enumerate
both primitive G-words and primitive R-words.  Specifically:

\begin{theorem}
\label{bijections}
Let $n\geq 2$.  Then:
\begin{enumerate}
\item The primitive G-words on $[n]$ are equinumerous with the decreasing 012-trees on vertex set $[n-2]$.
\item The primitive R-words on $[n]$ are equinumerous with the decreasing 012-trees on vertex set $[n-2]$.
\end{enumerate}
\end{theorem}

To prove this theorem, we construct explicit bijections between
G-words and decreasing 012-trees (Theorem~\ref{G-bijection}) and
between R-words and decreasing 012-trees (Theorem~\ref{R-bijection}).
Our constructions are of the same ilk as Donaghey's bijection \cite{Donaghey}
between decreasing 012-trees on $[n]$ and \emph{updown permutations},
i.e., permutations $w=w_1w_2\cdots w_n\in\Sym_n$ such that $w_1<w_2>w_3<\cdots$.

Together with \eqref{num-gen}, Theorem~\ref{bijections} enumerates the generators of
the graded-lex and reverse-lex initial ideals of $I_n$ in each degree; for instance,
$I_6$ is generated by $\binom{6}{4}\cdot1=15$ cubic monomials,
$\binom{6}{5}\cdot2=12$ quartics, and $\binom{6}{6}\cdot5=5$ quintics.

We start with three lemmas describing the recursive structure of G- and R-words.

\begin{lemma}
\label{splittable}
Let $n\geq 3$, let $w\in\Sym_n$, and let $k=w_{n-2}^{-1}$.  Define
words $w_L,w_R$ by
  $$
  w_L = w_1 w_{k-1} w_{k-2}\cdots w_3 w_2 w_k, \qquad
  w_R = w_n w_{k+1} w_{k+2}\cdots w_{n-2} w_{n-1} w_k.
  $$
Then:
\begin{enumerate}
\item If $w$ is a primitive G-word, then so are $w_L$ and $w_R$.
\item If $w$ is a primitive R-word, then so are $w_L$ and $w_R$.
\end{enumerate}
\end{lemma}

\begin{proof}
We will show that if $w$ is a primitive G-word, then so is $w_L$; the other
cases are all analogous.
If $n=3$, then the conclusion is trivial.  Otherwise, we have
$2\leq k\leq n-2$ by definition of a G-word.
If $k=2$, then $w_L=w_1w_2$, while if $k=3$, then $w_L=w_1w_3w_2$;
in both cases the conclusion follows by inspection.
Now suppose that $k\geq 4$.  Then the definition of~$k$ implies that $w_L$ satisfies \textbf{(G1)},
and if $w_{k-1}<w_2$ then $w[1,k]$ is a G-word, contradicting
the assumption that $w$ is a primitive G-word.  Therefore $w_L$ is a G-word.
Moreover, $w_L[i,j]\approx w[k+1-j,k+1-i]^*$ for every $[i,j]\subsetneq[k]$.
No such subword of $w$ is a G-word, so $w_L$ is a primitive G-word as desired.
\end{proof}

\begin{lemma}
\label{type-one} Let $n\geq 3$ and $x\in\Sym_{n-1}$.
\begin{enumerate}
\item If $x$ is a primitive G-word, then so is
  $$w = (n,\;n-2,\;x_2,\;x_3,\;\dots,\;x_{n-2},\;n-1).$$
\item If $x$ is a primitive R-word, then so is
  $$w = (n,\;x_{n-2},\;x_{n-3},\;\dots,\;x_2,\;n-2,\;n-1).$$
\end{enumerate}
\end{lemma}

\begin{proof}
Suppose that $x$ is a primitive G-word.
By construction, $w$ is a G-word in $\Sym_n$.  Let $w[i,j]$ be any proper subword of~$w$.
Then:
\begin{itemize}
  \item If $i\geq 3$, or if $i=2$ and $j<n$, then $w[i,j]=x[i-1,j-1]$ is not a G-word.
  \item If $i=2$ and $j=n$, then $w_i<w_j$ but $w_{i+1}=x_2>w_{j-1}=x_{n-2}$
(because $x$ is a G-word), so $w[i,j]$ is not a G-word.
  \item If $i=1$, then $j<n$, but then $w_{i+1}\geq w_j$, so~$w[i,j]$ is not a G-word.
\end{itemize}
Therefore $w$ is a primitive G-word.
The proof of assertion~(2) is similar.
\end{proof}

\begin{lemma}
\label{type-two}
Let $n\geq 4$, and let $P,Q$ be subsets of $[n]$ such that
  $$p=|P|\geq 3,\quad q=|Q|\geq 3,\quad P\cup Q=[n], \quad\text{and}\quad P\cap Q=\{n-2\}.$$
Let $x\in\Sym_P$ and $y\in\Sym_Q$ such that $x_p=n-2=y_q$ and $x_{p-1}>y_{q-1}$.
Then:
\begin{enumerate}
\item If $x$ and $y$ are primitive G-words, then so is
  $$w = (n, \; x_{p-1},\;\dots,\;x_2, \; n-2, \; y_2,\;\dots,\;y_{q-1}, \; n-1).$$
\item If $x$ and $y$ are primitive R-words, then so is
  $$w = (n, \; y_{q-1},\;\dots,\;y_2, \; n-2, \; x_2,\;\dots,\;x_{p-1}, \; n-1).$$
\end{enumerate}
\end{lemma}

\begin{proof}
Suppose that $x$ and $y$ are primitive G-words.
By construction, $w$ is a G-word.  We will show that
no proper subword $w[i,j]$ of $w$ is a G-word.  Indeed:
\begin{itemize}
  \item If $i<p<j$, then $w[i,j]$ cannot satisfy \textbf{(G1)}.
  \item If $i\geq p$, then either $[i,j]=[p,n]$, when $w_i=n-2<w_j=n-1$
and $w_{i+1}=y_2\geq w_{j-1}=y_{q-1}$ (because $y$ is a G-word),
or else $[i,j]\subsetneq[p,n]$, when $w[i,j]\approx y[i-p+1,j-p+1]$.  In either case,
$w[i,j]$ is not a G-word.
  \item Similarly, if $j\leq p$, then either $[i,j]=[1,p]$, when $w_i>w_j$ and
$w_{i+1}=x_{p-1}\leq w_{j-1}=x_2$ (because $x$ is a G-word),
or else $[i,j]\subsetneq[1,p]$, when $w[i,j]^* \approx x[p-j+1,p-i+1]$.  In either case,
$w[i,j]$ is not a G-word.
\end{itemize}
Therefore, $w$ is a primitive G-word.
The proof of assertion~(2) is similar.
\end{proof}

We pause to point out an elementary property about primitive G-words, which
is not necessary for the sequel, but is easy to observe from \eqref{small-examples}
and can be proved by an argument similar to the preceding lemmas.

\begin{proposition}
\label{penultimate}
Let $n\geq 2$ and let $w\in\MGW_n$.  Then $w_{n-1}=1$.
\end{proposition}

\begin{proof}
For $n\leq 6$, the result follows by inspection from \eqref{small-examples}.
Otherwise, let $i=w_1^{-1}$.  Note that $i\not\in\{1,2,n\}$
by the definition of G-word.
Suppose that $i\neq n-1$ as well.
By replacing $w$ with $w^*$ if necessary, we may assume that
$w_{i-1}<w_{i+1}$.  Let
$A=\{j\in[1,i-2]\st w_j>w_{i+1}\}$.  In particular
$\{1\}\subseteq A\subseteq[1,i-2]$.  Let
$k=\max(A)$.  Then
  \begin{align*}
  w_{k}   &= \max\{w_{k}, w_{k+1}, \dots, w_{i+1}\},\\
  w_{i+1} &= \max\{w_{k+1}, \dots, w_{i+1}\},\\
  w_{k+1} &> w_i=1.
  \end{align*}
So $w[k,i+1]$ is a G-word.  It is a proper
subword of $w$ because $i+1\leq n-1$, and its length is $i+2-k\geq i+2-(i-2)=4$.
Therefore $w\not\in\MGW_n$.
\end{proof}

For the rest of the paper, let $P$ be a finite subset of $\Nn$, let
$n=|P|$, and let $m=\max(P)$.  Define
  \begin{align*}
  \MGW'_P &= \{w\in\Sym_P \st (m+2,w,m+1)\in\MGW\},\\
  \MRW'_P &= \{w\in\Sym_P \st (m+2,w,m+1)\in\MRW\}.
  \end{align*}
The elements of $\MGW'_P$ (resp.,\ $\MRW'_P$) should be regarded as primitive G-words
(resp.,\ primitive R-words) on $P\cup\{m+1,m+2\}$, from which the first and
last digits have been removed.

We now construct a bijection between $\MGW'_P$
and the decreasing 012-trees $\DD_n$ on vertex set $[n]$.  If $P=\0$, then both
these sets trivially have cardinality 1, so we assume henceforth that $P\neq\0$.
Since the cardinalities of $\MGW'_P$ and $\DD_P$ depend only on $|P|$,
this theorem is equivalent to the statement that the primitive G-words on $[n]$
are equinumerous with the decreasing 012-trees on vertex set $[n-2]$, which
is the first assertion of Theorem~\ref{bijections}.

Let $w\in\MGW'_P$ and $k=w_m^{-1}$.
Note that if $n>1$, then $w_n<w_1\leq m$, so $k\neq n$.
Define a decreasing 012-tree $\phi_G(w)$ recursively (using the notation of
Definition~\ref{decreasing-tree}) by
  \begin{equation*}
  \phi_G(w) = \begin{cases}
    [m] & \text{ if } n=1;\\
    [m,\phi_G(w[2,n])] & \text{ if } n>1 \text{ and } k=1;\\
    [m,\phi_G(w[1,k-1]^*),\phi_G(w[k+1,n])] & \text{ if } n>1 \text{ and } 2\leq k\leq n-1.
  \end{cases}
  \end{equation*}

Now, given $T\in\DD_P$, recursively define a word
$\psi_G(T)\in\Sym_P$ as follows.
\begin{itemize}
\item If $T$ consists of a single vertex $v$, then $\psi_G(T)=m$.
\item If $T=[m,T']$, then $\psi_G(T)=(m,\psi_G(T'))$.
\item If $T=[m,T',T'']$ with $\min(P)\in T''$, then $\psi_G(T)=(\psi_G(T')^*,m,\psi_G(T''))$.
\end{itemize}

For example, let $T$ be the decreasing 012-tree shown in Definition~\ref{decreasing-tree}.
Then
\begin{align*}
\psi_G(T) &= \psi_G\left([6,[5,[4],[2]],[3,[1]]]\right)\\
&= \left( \psi_G([5,[4],[2]])^*,\ 6,\ \psi_G([3,[1]]) \right)\\
&= \left( (452)^*,\ 6,\ 31 \right)\\
&= 254631
\end{align*}
which is an element of $\MGW_6$ because, as one may verify, 82546317 is a primitive G-word.
Meanwhile, $\phi_G(254631)=T$.

\begin{theorem}
\label{G-bijection}
The functions $\phi_G$ and $\psi_G$ are bijections $\MGW'_n \to \DD_n$
and $\DD_n \to \MGW'_n$ respectively.
\end{theorem}

\begin{proof}
First, we show by induction on~$n=|P|$ that $\psi_G(T)\in\MGW'_P$.
This is clear if $n=1$; assume that it is true for all decreasing 012-trees
on fewer than~$n$ vertices.  If $T=[m,T']$, then $\psi_G(T)\in\MGW'_P$ by Lemma~\ref{type-one},
and if $T=[m,T',T'']$, then $\psi_G(T)\in\MGW'_P$ by Lemma~\ref{type-two}.

Next, we show that $\phi_G$ and $\psi_G$ are mutual inverses.
Again, we proceed by induction on~$n$.
The base case $n=1$ is clear, so we assume henceforth $n>1$.

Suppose inductively that $\psi_G(\phi_G(x))=x$ for all $x$ with $|x|<|w|$.
Suppose $w_k=m=\max(P)$.  If $k=1$, then by induction
  $$\psi_G(\phi_G(w)) = \psi_G([m,\phi_G(w[2,n])])
  = (m,\psi_G\phi_G(w[2,n]))
  = (m,w[2,n]) = w$$
while if $2\leq k\leq n-1$, then
  \begin{align*}
  \psi_G(\phi_G(w)) &= \psi_G([m,\phi_G(w[1,k-1]^*),\phi_G(w[k+1,n])])\\
  &= (\psi_G(\phi_G(w[1,k-1]^*))^*, m, \psi_G(\phi_G(w[k+1,n])))\\
  &= ((w[1,k-1]^*)^*, m, w[k+1,n]) = w.
  \end{align*}

On the other hand, suppose inductively that $\phi_G(\psi_G(U))=U$ for every
tree $U\in\DD_P$ with $|U|<n$.
If $T=[m,T']$, then
  $$\phi_G(\psi_G(T))
  = \phi_G(m,(\psi_G(T')))
  = [m,\phi_G(\psi_G(T'))]
  = [m,T'] = T$$
while if $T=[m,T',T'']$ with $\min(P)\in T''$, then
  \begin{align*}
  \phi_G(\psi_G(T))
  &= \phi_G((\psi_G(T')^*,m,\psi_G(T'')))
  = [m,\phi_G(\psi_G(T')),\phi_G(\psi_G(T''))]\\
  &= [m,T',T''] = T
  \end{align*}
as desired.
\end{proof}

Next, we construct the analogous bijections for primitive R-words.
Let $w\in\MRW'_P$ with $k=w^{-1}_m$.  Note that if $n>1$, then
$w_1<w_n\leq m$, so $k\neq 1$.  Define
a decreasing 012-tree $\phi_R(w)$ recursively by
  \begin{equation*}
  \phi_R(w) = \begin{cases}
    [m] & \text{ if } n=1;\\
    [m,\phi_R(w[1,n-1]^*)] & \text{ if } n>1 \text{ and } k=n;\\
    [m,\phi_R(w[1,k-1]^*),\phi_R(w[k+1,n])] & \text{ if } n>1 \text{ and } 2\leq k\leq n-1.
  \end{cases}
  \end{equation*}

Now, given $T\in\DD_P$, we recursively define a word
$\psi_R(T)\in\Sym_P$ as follows.
\begin{itemize}
\item If $T$ consists of a single vertex $v$, then $\psi_R(T)=v$.
\item If $T=[v,T']$, then $\psi_R(T)=(\psi_R(T')^*,v)$.
\item If $T=[v,T',T'']$, and the last digit of $\psi_R(T')$ is less than
the last digit of $\psi_R(T'')$, then $\psi_R(T)=(\psi_R(T')^*,v,\psi_R(T''))$.
\end{itemize}

Again, if $T$ is the decreasing 012-tree shown in Definition~\ref{decreasing-tree},
then
\begin{align*}
\psi_R(T) &= \psi_R\left([6,[3,[1]],[5,[4],[2]]]\right)\\
&= \left( \psi_R([3,[1]])^*,\ 6,\ \psi_R([5,[2],[4]]) \right)\\
&= \left( (13)^*,\ 6,\ 254 \right)\\
&= 316254
\end{align*}
which is an element of $\MRW_6$ because, as one may verify, 83162547 is a primitive R-word.
Meanwhile, $\phi_R(316254)=T$.

\begin{theorem}
\label{R-bijection}
The functions $\phi_R$ and $\psi_R$ are bijections $\MRW'_n \to \DD_n$
and $\DD_n \to \MRW'_n$ respectively.
\end{theorem}

\begin{proof}
First, we show by induction on~$n=|P|$ that $\psi_R(T)\in\MRW'_P$.
This is clear if $n=1$, so assume that it is true for all decreasing 012-trees
on fewer than~$n$ vertices.  If $T=[v,T']$, then $\psi_R(T)\in\MRW'_P$ by Lemma~\ref{type-one},
and if $T=[v,T',T'']$, then $\psi_R(T)\in\MRW'_P$ by Lemma~\ref{type-two}.

We have now constructed functions
  $$
  \phi_R: \ \MRW'_n \to \DD_n,\qquad
  \psi_R: \ \DD_n \to \MRW'_n.
  $$
It remains to show that they are mutual inverses, which we do by induction on~$n$.
The base case $n=1$ is clear, so we assume henceforth $n>1$.

Suppose inductively that $\psi_R(\phi_R(x))=x$ for all $x$ with $|x|<|w|$.
Suppose $w_k=m=\max(P)$.  If $k=n$, then by induction
  \begin{align*}
  \psi_R(\phi_R(w)) &= \psi_R([m,\phi_R(w[1,n-1]^*)])
  = (\psi_R(\phi_R(w[1,n-1])),m)\\
  &= (w[1,n-1],m) = w
  \end{align*}
while if $2\leq k\leq n-1$, then
  \begin{align*}
  \psi_R(\phi_R(w)) &= \psi_R([m,\phi_R(w[1,k-1]^*),\phi_R(w[k+1,n])])\\
  &= (\psi_R(\phi_R(w[1,k-1]^*))^*, m, \psi_R(\phi_R(w[k+1,n])))\\
  &= ((w[1,k-1]^*)^*, m, w[k+1,n]) = w.
  \end{align*}

On the other hand, suppose inductively that $\phi_R(\psi_R(U))=U$ for all $U\in\DD_P$ with $|P|<n$.
If $T=[m,T']$, then
  \begin{align*}
  \phi_R(\psi_R(T))
  &= \phi_R((\psi_R(T')^*,m))
  = [m,\phi_R((\psi_R(T')^*)^*)]
  = [m,\phi_R(\psi_R(T'))]\\
  &= [m,T'] = T
  \end{align*}
while if $T=[m,T',T'']$ with $\psi_R(T')_1<\psi_R(T'')_1$, then
  \begin{align*}
  \phi_R(\psi_R(T))
  &= \phi_R((\psi_R(T')^*,m,\psi_R(T'')))
  = [m,\phi_R(\psi_R(T')),\phi_R(\psi_R(T''))]\\
  &= [m,T',T''] = T
  \end{align*}
as desired.
\end{proof}

\bibliographystyle{amsalpha}

\end{document}